\documentclass{preprint} 

\usepackage[overload]{abraces} 		
\usepackage{enumerate} 				
\usepackage{amsthm}
\usepackage{stmaryrd} 
\usepackage{cleveref}
\usepackage{grundgesetze}				    
\usepackage{bussproofs}
\usepackage{multicol}

\allowdisplaybreaks
\usepackage{tikz-cd}
\usetikzlibrary{patterns}
\tikzcdset{every label/.append style = {font = \small}} 
\usetikzlibrary{decorations.markings}

\newtheorem{theorem}{Theorem}[section]
\newtheorem{lemma}[theorem]{Lemma}
\newtheorem{corollary}[theorem]{Corollary}
\newtheorem{definition}[theorem]{Definition}

\usepackage{setspace}

\newcommand{\bool}{\{0_D,1_D\}}
\newcommand{\M}{\mathcal{M}}
\newcommand{\var}{\mathsf{var}}
\newcommand{\term}{\mathsf{term}}
\newcommand{\form}{\mathsf{form}}

\renewcommand{\nat}{\mathbb{N}}

\newcommand{\GGA}{\textsf{GGA1}}

\newcommand{\ia}{A1}
\newcommand{\iia}{A2}
\newcommand{\iiia}{A3a}
\newcommand{\iiie}{A3e}
\newcommand{\iv}{A4}

\newcommand{\ahor}{A5}
\newcommand{\aneg}{A6}
\newcommand{\aimp}{A7}
\newcommand{\aquant}{A8}
\newcommand{\aid}{A9}

\newcommand{\hypsyll}{A10}
\newcommand{\dilemma}{A11}
\newcommand{\fug}{A12}

\renewcommand{\mp}{MP}
\newcommand{\ug}{UG}

\DeclareMathSymbol{\supvar}{3}{matha}{"81} 

\renewcommand{\to}{\supvar}

\let\phi\varphi

\makeatletter
\newcommand\mytag[2][]{%
	\def\@currentlabel{#2}%
	(#2)\label{#1} 
}
\makeatother

\title{Frege's logic in {\gga} is nonclassical}

\journal{Notre Dame Journal of Formal Logic}

\volume{67}

\issue{(3), }

\pages{pp. 437--455}

\when{2026}

\author{Bruno Bentzen}

\doi{10.1215/00294527-2026-0008}

\department{School of Philosophy}

\university{Zhejiang University}

\country{Hangzhou, China}

\email{bbentzen@zju.edu.cn}

\begin{document}
	
\maketitle

\begin{abstract}
	I show that the axioms and rules of inference about the horizontal in the subsystem of \gga~without value-ranges fail to capture its intended interpretation. That is, I build a nonstandard model in which $-a$ is the true iff $a$ is not the false, thereby confirming a conjecture made by Landini that $a \to (-a) = a$ is not provable in this subsystem. This motivates our main result, namely, soundness and completeness proofs for the first-order fragment without value-ranges with respect to a new nonclassical semantics of boolean truncation. Simply put, the widespread assumption that~\gga's first-order logic is classical first-order logic is mistaken.
\end{abstract}


\nocite{frege1962grundgesetze}

\section{Introduction} \label{intro}

Recent studies have questioned the received view that the role of the horizontal in the formalism of \textit{Grundgesetze der Arithmetik} is superfluous.\footnote{
	This received view can be traced back to~\citet[p.~315]{dummett1973fregelanguage}. I return to it in~\Cref{conclusion}. } 
\cite{landini2022stipulations} conjectures that $a \to (-a)=a$ is independent relative to the system without value-ranges. 
This is true in the intended interpretation of the horizontal, where $-a$ is the true iff $a$ is the true. But Frege's axiomatization seems to say too little about the horizontal to decide this statement. 
Clearly, the horizontal is indispensable if this statement would not qualify as a theorem. 
Independently, \cite{bellucci2023frege} study the amalgamation of horizontals, arguing that the horizontal is a crucial element of Frege's system because of the way it contributes to the sense of the sentences in which it is included. 
But as proof-theoretic irredundance should not~amount to questions of sense, their claim is better understood as saying that the horizontal is semantically irredundant with respect to its intended interpretation. What if the intended~interpretation of the horizontal cannot be captured by the formalism? 

\newpage

My first aim in this paper is to show that the axioms and rules of inference governing the horizontal in the subsystem of \gga~without value-ranges admit a nonstandard interpretation in which $-a$ is the true iff $a$ is not the false, thereby confirming Landini's conjecture about the unprovability of $a \to (-a)=a$. My second and primary aim is to answer the question of what semantics its first-order fragment is complete with respect to. I isolate a first-order system called~$\GGA$, designed to make explicit all the axioms and rules that govern the horizontal in \gga, and prove soundness and completeness in terms of a new kind of nonclassical semantics based on boolean truncation. 
This suggests that the widespread assumption, well~illustrated in~\cite{parsons1987consistency}, that the first-order logic of \gga~amounts to classical first-order logic is mistaken. 
The remainder of this paper consists of three main sections. 
The independence proof is given in \Cref{sec:independence} and lays the~foundation for our study of \GGA~in \Cref{sec:completeness}. 
I provide some concluding remarks in \Cref{conclusion} and comment on the philosophical relevance of our result.

\section{The independence proof} \label{sec:independence}

\cite{landini2022stipulations} makes the intriguing observation that no axiom is given to accommodate the semantic stipulation that the negation of the number~$2$ is the true in~\gga~\S6. 
He argues, based on some pieces of textual evidence, that the first step toward filling this gap in the formalism is realizing that Frege informally held the following:\footnote{
	See \citet[p.~352, \S4]{landini2022stipulations}. I postpone a discussion of his argument to \Cref{conclusion}.} 
\begin{equation}\tag{A} \label{A}
	a \to (-a)=a.
\end{equation}
Simply put, this says that if $a$ is true then the horizontal of $a$ is identical to $a$. 
Surely this would be a natural theorem to have given the intended interpretation of the horizontal. 
But upon closer examination we see that (\ref{A}) is not proved anywhere in {\gga} and Frege never explicitly states it as an axiom of the system. 
Landini sees ``very strong indications of its independence'' relative to the system of \gga~minus Basic Law V, 
because of the way horizontals work in Frege's formalism, but argues, in an interesting appendix, that ``modern techniques for proving independence cannot apply in the context of Frege's horizontal and Amalgamation''~\cite[p.~352, fn.~8]{landini2022stipulations}.  
For him, the presence of the horizontal is what poses unique new issues which have no known solution. He argues that if we simply follow the usual independence proof methods, we will ignore the presence of the horizontal and its amalgamation in the system~\cite[pp.~379--381]{landini2022stipulations}. 
But this need not be the case. In the remainder of this section, I show that an independence proof for (\ref{A}) can be given if we explicitly consider the horizontal as a logical operator with a nonstandard semantics.

To begin with, we must note that, unlike the formalism of \bs, where the logical connectives flank formulas, \gga~is based on a term language, where the logical connectives are term-forming functions returning objects. 
This is a point recently stressed by~\cite{landini2012frege}, \cite{cadetpanza2015logical}, and \cite{bentzen2023fregett}, but missed in most modern reconstructions of Frege's system. 
In the first two works, a turnstile notation is used to express the formula resulting from a term $a$ as $\vdash a$. 
I will not follow this convention here, since it invites a dangerous confusion between judgments and propositions.\footnote{
	This has been stressed elsewhere~\citep[\S\S4.2--4.3]{bentzen2023fregett}. See also~\cite{kim2020horizontal}.} 
I will turn to a more rigorous treatment of the language of \gga~in~\Cref{language}. 
For now, 
we can speak of formulas and terms interchangeably if we stipulate that the truth value of a formula is determined by that of its corresponding term. It will suffice to bear in mind that any independence proof would have to assign truth values to terms and not formulas.

\subsection{Propositional fragment (I and IV)} \label{axioms}

We begin with the propositional fragment of~\gga. Strictly speaking, it is not reasonable to speak of such a fragment because we are always dealing with terms. But the division will serve a methodological purpose. 
Under this group we may consider as connectives the horizontal, negation, implication, and identity.  We include as axioms Basic Laws~I and IV and, as rules of inference, the permutation of subcomponents, amalgamation of identical subcomponents, amalgamation of horizontals, and generalized modus ponens, hypothetical syllogism, contraposition, and dilemma, in Cook's~(\citeyear{sep-frege-logic}, \S3.4) terminology.

Since (\ref{A}) is true under the intended interpretation of the horizontal, we must appeal to a nonstandard interpretation in order to falsify it. 
Indeed, the worry expressed by Landini concerning the availability of independence proofs only manifests itself because he is explicitly adhering to the intended interpretation of the horizontal: 

\begin{quote}
	Where $x$ is the designated value, $val(`\vdash (−a)') = x$, if $val(`\vdash a') = x$; and for every other object $\Delta$ such that $val(`\vdash a’) = \Delta$, we are to have $val(`\vdash (−a)’) = y$, for some fixed $y$. [...] \emph{If the designated
	value is to be $0$ (in our example above) then the horizontal should return some fixed $y$ for any input
	other than $0$, and return $0$ for an input of $0$. If we choose $y$ as 1, then in respect of Amalgamation and
	the horizontal the table will alter any subordinate $2$ in the original table to a $1$.} This, I fear, spoils all
	hope of applying the usual techniques for independence proofs. \cite[p.~381, emphasis ours]{landini2022stipulations} 
\end{quote}

\paragraph{Logical connectives} 

The basic intuition behind the truth tables below is that two distinct truth-values play the role of the true, namely, $1$ and $2$. Basically, the horizontal truncates truth by conflating $1$ and $2$ into $1$. Only in this admittedly weak sense is the standard interpretation captured, given that both $1$ and $2$ are the true, but $1$ always has precedence over $2$: 


\begin{center}
	\begin{tabular}{c | c}
		$-$ &  \\
		\hline
		0 & 0  \\
		1 & 1 \\
		2 & 1 \\ 
	\end{tabular}
	\quad
	\begin{tabular}{c | c}
		$\neg$ &  \\
		\hline
		0 & 1  \\
		1 & 0 \\
		2 & 0 \\ 
	\end{tabular}
	\quad
	\begin{tabular}{c | c  c c}
		$=$ & 0 & 1 & 2 \\
		\hline
		0 & 1 & 0 & 0 \\
		1 & 0 & 1 & 0 \\
		2 & 0 & 0 & 1 \\ 
	\end{tabular}
	\quad
	\begin{tabular}{c | c  c c}
		$\to$ & 0 & 1 & 2 \\
		\hline
		0 & 1 & 1 & 1 \\
		1 & 0 & 1 & 1 \\
		2 & 0 & 1 & 1 \\ 
	\end{tabular}
\end{center}


From a formal perspective, the above tables for the horizontal, negation, and implication capture precisely the alternative versions contemplated by~\citet[\S4]{berg2017propositional}, except that they think of the horizontal as true iff the argument is nonfalse. 
They observe that this alternative semantics for the horizontal and negation can be understood in light of Frege's trans-sortal identification of truth values with value-ranges in \gga~\S10, but maintain that they should not be regarded as ``an actual \textit{alternative} to the primitive propositional connectives used by Frege in constructing the theorems of \textit{Grundgesetze}''. 
Our independence proof shows that their claim is too modest: as far as the axioms and rules of the formalism are concerned, these alternative interpretations are just as good as the original ones.

\paragraph{Axioms and rules of inference} 
Our designated value is $1$. We thus need to show that all axioms have value~$1$ and that all rules of inference preserve this designated value. 
Basic Law I is in fact a group of two axioms $a \to (b \to a )$ and $a \to a$, while Basic Law IV states that $\neg ((-a)=(\neg b)) \to ((-a) = (-b))$. 
Their validation is easy to see in the tables below:

\begin{center}
	\begin{tabular}{c c c c c }
		$a $ & $\to$ & $(b$ & $\to$ & $a )$ \\
		\hline
		0 & 1 & 	0 	& 1 & 0\\
		0 & 1 & 	1 	& 0 & 0\\
		0 & 1 & 	2 	& 0 & 0\\
		1 & 1 & 	0 	& 1 & 1\\
		1 & 1 & 	1 	& 1 & 1\\
		1 & 1 & 	2 	& 1 & 1\\
		2 & 1 & 	0 	& 1 & 2\\
		2 & 1 & 	1 	& 1 & 2\\
		2 & 1 & 	2 	& 1 & 2\\
	\end{tabular}
	\qquad 
	\begin{tabular}{c c c  }
		$a $ & $\to$ & $a $ \\
		\hline
		0 & 1 &  0\\
		1 & 1 &  1\\
		2 & 1 &  2\\
	\end{tabular}
\end{center}
\begin{center}
	\begin{tabular}{c c c c c c c c c c c c}
		$(\neg $ & $(-$ & $a)$ & $=$ & $(\neg$ & $b))$ & $\to$ & $((-$ & $a )$ & = & $(-$ & $b ))$ \\
		\hline
		1 & 0 & 	0 	& 0 & 1& 	0& 	1 & 0 & 	0& 		1& 0& 	0\\
		0 & 0 & 	0 	& 1 & 0& 	1& 	1 & 0 & 	0& 		0& 1& 	1\\
		0 & 0 & 	0 	& 1 & 0& 	2& 	1 & 0 & 	0& 		0& 1& 	2\\
		0 & 1 & 	1 	& 1 & 1& 	0& 	1 & 1 & 	1& 		0& 0& 	0\\
		1 & 1 & 	1 	& 0 & 0& 	1& 	1 & 1 & 	1& 		1& 1& 	1\\
		1 & 1 & 	1 	& 0 & 0& 	2& 	1 & 1 & 	1& 		1& 1& 	2\\
		0 & 1 & 	2 	& 1 & 1& 	0& 	1 & 1 & 	2& 		0& 0& 	0\\
		1 & 1 & 	2 	& 0 & 0& 	1& 	1 & 1 & 	2& 		1& 1& 	1\\
		1 & 1 & 	2 	& 0 & 0& 	2& 	1 & 1 & 	2& 		1& 1& 	2\\
	\end{tabular}
\end{center}

Since our truth tables of the logical connectives agree with classical semantics, 
we see that all rules of inference not directly involving the horizontal, such as 
permutation of subcomponents, amalgamation of identical subcomponents, and generalized modus ponens, hypothetical syllogism, contraposition, and dilemma, 
are designated-value preserving. 
So we only have to check that all rules of amalgamation of horizontals are respected. The rules are summarized in \gga~\S48 but introduced much earlier. 
In~\S6, Frege introduces the transition from $- (- a)$ to $-a$ and writes that $\neg a$ refers the same as $\neg(-a)$, $- (\neg a)$, and $- \neg (- a)$. Then,~in~\S12, he adds that $a \to b$ is the same as $- ((-a) \to (-b))$. Frege does not list an amalgamation rule for identity in \S48, but in $\S50$ he derives theorem~IIIi that states that $(a = b) = -(a = b)$ using second-order quantification and Basic Law III. 
The amalgamations are expressed in the tables below in terms of logical equivalences:

\begin{center}
	\begin{tabular}{c c | c c c }
		$-$ & $a$ & $-$ & $(-$ & $a )$ \\
		\hline
		0 & 0 & 0 & 0 & 0\;\\
		1 & 1 & 1 & 1 & 1\;\\
		1 & 2 & 1 & 1 & 2\; \\
	\end{tabular}
	\qquad 
	\begin{tabular}{c c | c c c c }
		$\neg$ & $a$ & $-$ & $(\neg$ & $(-$ & $a ))$ \\
		\hline
		1 & 0\; & 1 & 1 & 0& 0\; \\
		0 & 1\; & 0 & 0 & 1& 1\; \\
		0 & 2\; & 0 & 0 & 1& 2\; \\
	\end{tabular}

	\begin{tabular}{c c c | c c c c }
		$a $ & $=$ & $b$ & $(-$ & $( a$ & $=$ & $ b))$ \\
		\hline
		0 & 1 & 	0 	& 1 & 0 & 1 & 0 \\
		0 & 0 & 	1 	& 0 & 0 & 0 & 1 \\
		0 & 0 & 	2 	& 0 & 0 & 0 & 2 \\
		1 & 0 & 	0 	& 0 & 1 & 0 & 0 \\
		1 & 1 & 	1 	& 1 & 1 & 1 & 1 \\
		1 & 0 & 	2 	& 0 & 1 & 0 & 2 \\
		2 & 0 & 	0 	& 0 & 2 & 0 & 0 \\
		2 & 0 & 	1 	& 0 & 2 & 0 & 1 \\
		2 & 1 & 	2 	& 1 & 2 & 1 & 2 \\
	\end{tabular}
	\qquad 
	\begin{tabular}{c c c | c c c c c c}
		$a $ & $\to$ & $b$ & $-$ & $(-$ & $a)$ & $\to$ & $(-$ & $b)$ \\
		\hline
		0 & 1 & 	0 	& 1 & 0 & 0\; & 1 & 0 & 0\; \\
		0 & 1 & 	1 	& 1 & 0 & 0\; & 1 & 1 & 1\; \\
		0 & 1 & 	2 	& 1 & 0 & 0\; & 1 & 1 & 2\; \\
		1 & 0 & 	0 	& 0 & 1 & 1\; & 0 & 0 & 0\; \\
		1 & 1 & 	1 	& 1 & 1 & 1\; & 1 & 1 & 1\; \\
		1 & 1 & 	2 	& 1 & 1 & 1\; & 1 & 1 & 2\; \\
		2 & 0 & 	0 	& 0 & 1 & 2\; & 0 & 0 & 0\; \\
		2 & 1 & 	1 	& 1 & 1 & 2\; & 1 & 1 & 1\; \\
		2 & 1 & 	2 	& 1 & 1 & 2\; & 1 & 1 & 2\; \\
	\end{tabular}
\end{center}

It is easy to check that if $a$ is a theorem, then $a$ has the same truth table as $-a$: that is because the horizontal is the identity on the designated value $1$.  

\paragraph{Landini's axiom (A)} As far as the propositional fragment is concerned, to complete the proof we just have to show that $a \to ((-a)=a)$ is not valid under our interpretation. It fails precisely because there are two nonidentical objects that both play the role of the true, thus there is a model where $a$ has the value $2$ but $(-a)$ is evaluated as $1$: 

\begin{center}
	\begin{tabular}{c c c c c c }
		$a$ & $\to$ & $((-$ & $a )$ & = & $a )$ \\
		\hline
		0 & 1 & 0 & 0& 1& 0\\
		1 & 1 & 1 & 1& 1& 1\\
		2 & 0 & 1 & 2& 0& 2\\
	\end{tabular}
\end{center}

\subsection{Quantificational fragment (IIa, IIb, III)} \label{axioms}

We extend the proof to the quantificational fragment that adds the first- and second-order universal quantifiers governed by Basic Laws IIa, IIb, and III as axioms, and concavity introduction and Roman letter elimination as rules of inference~\cite[\S\S3.4.5--6]{sep-frege-logic}. 
I omit not just Basic Law V but also VI as well, since its purpose is to handle the application of the definite description operator to value-ranges, which are absent from our consideration. 

\paragraph{Quantifers}

We shall write $a^\mu_{\mathcal{M}}$ to indicate the interpretation of the term~$a$ in the model~$\mathcal{M}$ under a certain variable assignment~$\mu$. The domain of our class of models is the set $\{0, 1, 2\}$ and thus models differ only in their interpretation of nonlogical function symbols. 
The truth tables given above in fact indicate the interpretation of terms without quantifiers. 
This makes it clear that the propositional fragment is inseparable from the quantificational one. 
The first- and second-order universal quantifiers are interpreted just as expected. Below we write $\mu, e/x$ for the variable assignment that assigns $e$ to $x$ but leaves all else unchanged:  
\begin{align*}
	(\forall x f(x))^\mu_{\mathcal{M}} &=
	\begin{cases*}
		1 & if $f(x)^{\mu, e/x}_{\mathcal{M}} \neq 0$ for all elements $e \in D$ \\
		0 & otherwise.
	\end{cases*} \\
	(\forall f M_x(f(x)))^\mu_{\mathcal{M}} &=
	\begin{cases*}
		1 & if $M_x(f(x))^{\mu, g/ f}_{\mathcal{M}} \neq 0$ for all functions $g \in D^D$ \\
		0 & otherwise.
	\end{cases*}
\end{align*}

\paragraph{Axioms and rules}

We write $\vDash a$ to mean that $(-a)^{\mu}_{\mathcal{M}}=1$ for all models $\mathcal{M}$ and variable assignments~$\mu$. We must show that Basic Laws IIa and IIb and Basic Law III are valid. 
Now, Basic Laws IIa and IIb state universal instantiation in conditional form for both universal quantifiers $\forall x f(x) \to f(a)$ and $\forall f M_x (f(x)) \to M_x (g(x))$. Clearly, both are valid, for it is never the case that the value of the antecedent is not $0$ but the consequent is $0$. 
Basic Law III states a more general version of Leibniz's principle of indiscernibility of identicals in the form of $g(a=b) \to g(\forall f(f(a) \to f(b)))$. This axiom is valid as well given that $(a=b)^{\mu}_{\mathcal{M}} = 1$ iff $(\forall f (f(a) \to f(b)))^{\mu}_{\mathcal{M}} = 1$ and we have a function $g^{\mu}_{\mathcal{M}} : D \rightarrow D$.

Both concavity
introduction and Roman letter elimination preserve our designated value $1$,
because they are essentially counterparts of classical universal introduction and elimination, and, just as in the propositional case, our semantics of the universal quantifier agrees with its classical semantics in all contexts that do not involve the horizontal. 
The amalgamation of horizontals for the first-order universal quantifier is found in \gga~\S8 with the stipulation that the term $\forall x f(x)$ refers the same as $\forall x (- f(x))$ and as $-(\forall x f(x))$. It suffices for us to observe that the identity $\forall x f(x)^\mu_{\mathcal{M}} = - (\forall x (- f(x)))^\mu_{\mathcal{M}}$ obtains. The same can be said of the second-order quantifiers Frege discusses in \gga~\S\S19--23.

\section{Completeness for the first-order fragment~\GGA} \label{sec:completeness}

So far I have been sloppy about the exact language, formal system, and semantics of the version of \gga~without value-ranges presupposed in our independence proof. But I hope I have said enough to motivate the pivotal question addressed in this paper: if the intended interpretation of the horizontal need not be fixed by the second-order formalism, then what semantics is its first-order fragment complete with respect to? 
Until now this question has not been seriously discussed in the literature, because it is simply taken for granted that the first-order fragment of {\gga} is classical first-order logic. 
For example,~the logical system considered by \cite{parsons1987consistency} in his famous proof of consistency of what he calls the ``first-order portion
of Frege's logical system'' is in fact nothing but classical~first-order logic with~identity augmented with Basic Law V. 

While Parsons's result is indeed technically interesting, 
it has very little connection to the actual first-order fragment of \gga~due to the omission of horizontals. It is not even meaningful~to~ask~whether (\ref{A}) is unprovable in Parsons's system. It features no axioms and rules involving the horizontal, although a logical connective for it is present in the language and it is easy to verify that (\ref{A}) is satisfied in his model because it adheres to the intended interpretation of the horizontal. The same applies to the results obtained by~\cite{heck1996consistency}, \cite{burgess1998consistent}, and many others building on Parsons's work.\footnote{
	Including \cite{wehmeier1999consistent,ferreira2002consistency, Ferreira2005amending, walsh2016frege}.}  
By~contrast, the horizontal is given the proper treatment it deserves in recent modern reconstructions proposed by~\cite{klement2002frege}, \cite{landini2012frege}, and \cite{cadetpanza2015logical}. 
The most relevant remarks on completeness I know are given by~Cadet and Panza. 
They prove that the propositional fragment of \gga~is deductively equivalent to classical propositional logic and reconstruct the second-order system without value-ranges, but assume the intended interpretation of the horizontal and do not isolate a first-order fragment. There has been no real motivation for a metatheoretical study of completeness for the first-order fragment of \gga~in the presence of the horizontal until now, when our independence proof of (\ref{A}) revealed that its intended interpretation is not forced by the formalism.

\subsection{The formal system~\GGA} \label{systemgga1}

In the rest of this paper, I propose the formal system \GGA~for the first-order fragment of \gga~without value-ranges and show that it is sound and complete with respect to a new class of algebras of boolean truncation. Unlike in the previous section, we will rigorously present its language, axiomatization, and semantics. Before we proceed, it should be mentioned that, as it has been observed multiple times in the literature, it is unlikely that Frege's logic in \gga~is fully amenable to modern reconstructions.\footnote{
	This point is nicely made by \citet{cadetpanza2015logical} and \citet{sep-frege-logic}. } 
There will always be design choices that need to be made to transpose his logic into the modern setting and there is no fully faithful embedding of his own two-dimensional notation. 
I do not see {\GGA} as an actual subsystem of Frege's logic in \gga, but as a modern reconstruction of its first-order fragment designed to accurately reflect the workings of the horizontal. 

Therefore, \GGA~is not intended to be a historically accurate system. I will not hesitate to deviate from Frege's design choices in \gga~if following him leads to unnecessary complications that divert attention from our study of the horizontal. 
For example, he includes redundant rules of inference to avoid lengthy derivations~\cite[p.~vi]{frege1962grundgesetze}. 
I will not follow this practice since I am mostly interested in reasoning about and not with~\GGA. 
In fact, to facilitate our metatheoretical analysis of the system, I will incorporate some of his rules of inference as axiom schemes whenever the same rules would remain admissible. 
Moreover, Frege allows a tacit use of amalgamation of horizontals, which I do not permit. 
Finally, I refrain from adopting Frege's two-dimensional notation because we are primarily interested in a metatheoretical study about a certain modern formal system, just as in \Cref{sec:independence}. The difference is that this formal system will now be fully spelled out for our completeness proof. As we will soon see, the nonstandard model from \Cref{sec:independence} also models \GGA.

\paragraph{Language} \label{language}

As already mentioned in \Cref{sec:independence}, \gga~features a term language. 
A signature $\sigma$ is a pair of a set of function symbols with an arity function. For simplicity, we may assume that the set of function symbols is countable and that each element has the form $f_n$ for $n \in \mathbb{N}$. We may write $f_n \in \sigma$ to mean that $f_n$ is a function symbol. We postulate a set~$\mathsf{var} = \{x_i \mid i \in \nat \}$ with an infinite stock of individual variables. We now have:

\begin{definition}
	Given any signature~$\sigma$, $\term(\sigma)$ is the smallest set closed under the conditions:
	\begin{enumerate}
		\item $x_i \in \term(\sigma)$, for all $x_i \in \mathsf{var}$;
		\item $f_i(a_1, ..., a_n) \in \term(\sigma)$, if $a_1,..,a_n \in \term(\sigma)$ and $f_i \in \sigma$ has arity $n$;
		\item $(a_1 = a_2) \in \term(\sigma)$, if $a_1,a_2 \in \term(\sigma)$;
		\item $(- a_1), (\neg a_1), (a_1 \to a_2) \in \term(\sigma) $, if $a_1, a_2 \in \term(\sigma)$;
		\item $(\forall x_i a_1) \in \term(\sigma)$, for all $x_i \in \var$, if $a_1 \in \term(\sigma)$.
	\end{enumerate}
\end{definition}

Our definition of terms follows~\citet[p.~46]{klement2002frege}, \citet[\S2.6]{landini2012frege} and \citet[\S\S3.2--3.3, 4.3]{cadetpanza2015logical} in all essential aspects. Notice also that we understand Frege's Roman object-letters as serving as individual variables, in agreement with \citet[\S1]{kemp1998propositions}, \citet[p.~38]{klement2002frege}, \citet[\S3]{heck2007frege}, and \citet[p.~32, \S5.1]{landini2012frege}. 
I stress that only terms obtained from clauses (1)--(2) are atomic. For better readability, we follow the usual bracket dropping conventions plus associativity to the right when there is no danger of ambiguity. 
We use the letters $a, b, ...$ as metavariables, noting that they never occur in the language. Instead of using an explicit substitution notation $f[a/x]$ we write $f(a)$, for $f, a \in \term(\sigma)$ and $x \in \var$, where it is always implicitly assumed that $a$ is free for $x$ in $f$. When $f$ is used without a subscript, it thus refers to a possibly open term and not a function symbol. 
Finally, it will be useful to adopt $\top$ as a short for~$\forall x (x=x)$ and $\bot$ as a short for $\neg \top$.

Theorems of modern first-order theories are stated using formulas and not terms. 
Since we wish \GGA~to stay as close to modern logic as possible, we still require a set of formulas to express our axioms and rules of inference. 
Of course, there will be only a single way of forming a formula by embedding a term into it. 
But technically speaking, we still need a convention to syntactically distinguish between terms and their induced formulas: 

\begin{definition}
	For any signature~$\sigma$, $\form(\sigma)$ is the smallest set closed under the condition:
	\begin{enumerate}
		\item $a' \in \form(\sigma)$, if $a \in \term(\sigma)$.
	\end{enumerate}
\end{definition}

That is, an apostrophe is used to indicate the injection of a term. 
Though this difference is essential, we shall omit the apostrophe to avoid a pedantic style.  Thus we do not follow~\cite{klement2002frege} in conflating terms and formulas into `well-formed expressions', but we do not always employ a notation to explicitly distinguish between them either, as \cite{landini2012frege} and \cite{cadetpanza2015logical} do (see~\Cref{sec:independence}). Readers should always be able to tell terms and formulas apart from the context. All our metavariables range over terms. 
From this point on, we operate on a fixed signature $\sigma$ for the sake of convenience.

\paragraph{Axiomatization} \label{subsec:axiomatization}  

To adapt the system of \gga~to a first-order setting, we remove Basic Law IIb and replace Basic Law III by adopting as primitive both theorems IIIa and IIIe that Frege derives from III in~\S50. This eliminates higher-order quantification. 
To avoid proof-theoretic redundancy we omit $a \to a$ from Basic Law~I. 
As a result, axiom schemes \ia--\iv~below incorporate Frege's formulas I, IIa, IIIa, IIIe, and IV. 
Here, to preserve the original numbering of \gga, we group IIIa and IIIe under the label A3, referring to them as \iiia~and \iiie~when there is danger of confusion. 
For convenience, we incorporate the rules of inference Frege lists in \gga~\S48 using simpler axioms or rules. 
That is, we employ \ahor--\aid~for the rules for amalgamation of horizontals plus IIIi and \hypsyll--\fug~and \mp--\ug~as a provably equivalent alternative to all other listed rules. 
We will see soon how to recover as admissible rules all Frege's original rules of inference overviewed in \Cref{sec:independence}.

\begin{definition}
	The first-order logic of \gga~without value-ranges~\GGA~contains exactly the axiom schemes and rules of inference below for $a, b, f \in \term(\sigma)$: 
	\begin{multicols}{2}
		\begin{enumerate}
			\item[\ia.] $ a \to (b \to a )$
			
			\item[\iia.] $\forall x f(x) \to f(a)$
			
			\item[A3a.]
			$a = b \to (f(a) \to f(b))$
			
			\item[A3e.]
			$a = a$ 
			
			\item[\iv.] $\neg ((-a)=(\neg b)) \to ((-a) = (-b))$
			
			\item[\ahor.] 
			$(-a) = (- (- a))$
			
			\item[\aneg.] 
			$(\neg a)
			=
			(-(\neg -a))$
			
			\item[\aimp.] 
			$(a \to b)
			=
			(- ((-a) \to (- b)))$
			
			\item[\aquant.] 
			$
			(\forall x f(x))
			=
			(- \forall x (- f(x)))
			$
			
			\item[\aid.] 
			$(a = b) =
			(-(a = b))
			$
			
			\item[\hypsyll.] $(a \to (b \to c))\to ((a \to b) \to (a \to c))$

			\item[\dilemma.] $(\neg b \to a) \to ((\neg b \to \neg a) \to b)$

			\item[\fug.] $\forall x (a \to f(x))\to (a \to \forall x f(x)))$

			\item[\mp.] \AxiomC{$a \to b$} 
			\AxiomC{$a$}
			\BinaryInfC{$b$} 
			\DisplayProof
			
			\item[\ug.] \AxiomC{$f(y)$} 
			\UnaryInfC{$\forall x f(x)$} 
			\DisplayProof 
		\end{enumerate}
	\end{multicols}
\end{definition}

Notice that the usual restrictions to \fug~and \ug~apply. In the former, we require that~$a$ lacks free occurrences of $x$ 
and, in the latter, that $y$ be a fresh variable free for $x$ in $f$. 
Also, the sequent $\Gamma \vdash a$ means that there exists a proof from premises in $\Gamma \subseteq \form(\sigma)$ and conclusion $a \in \form (\sigma)$, where proofs are understood as in the usual axiomatic way. 
Thus, all axioms and rules of inference considered here can be used under assumptions. This concludes~our axiomatization of the formal system~\GGA. However, 
a few general remarks for Frege scholars are in order to justify some of our design choices as well as the claim that \GGA~captures the first-order fragment of \gga~without value-ranges.

To begin with, Frege does not have hypothetical proofs in \gga. 
All premises used in a derivation must be found in the judgments involved, but his account of judgment as the acknowledgment of the truth of a thought from \S5 is strictly categorical, in the Kantian sense of a judgment not made under assumptions~\cite[\S\S4--5]{bentzen2023fregett}. Still, given that the deduction theorem is provable in \GGA, hypothetical proofs of~$\Gamma \vdash a$ can be useful to restore the original rules of inference from \gga~left out in our axiomatization. 
The proof of the deduction theorem is straightforward, given that the horizontal-identity free fragment of \GGA~consisting of \ia--\iia, \hypsyll--\fug, and \mp--\ug~is proof-theoretically equivalent to a standard Mendelson-style first-order axiomatic system without identity. But because the derivation takes place within a yet unfamiliar term-based formal system, it might be helpful to convince readers that we go through all basic steps in the same way.

\begin{theorem}[Deduction]
	If $\Gamma, a \vdash b$ then $\Gamma \vdash a \to b$.
\end{theorem}
\begin{proof}
	By induction on the structure of the derivation. As usual, the base case is covered by {\ia} and \mp, and, in the inductive case, \hypsyll~is employed for \mp, and \fug~for \ug. 
\end{proof}

Thus the fact that the logical operators connect terms rather than formulas in the deductive system does not introduce any additional proof-theoretic difficulties. (Notice that weakening follows for {\GGA} using the familiar argument too.)

\newpage
Using the deduction theorem, we can recover as admissible rules all the rules of inference listed by Frege in \gga~\S48 just as in classical logic. 
Since contraction and exchange hold in \GGA~because the premises of a sequent are sets, this yields Frege's rule of permutation of subcomponents, from $\vdash a_1 \to ...\ a_n \to a \to b$ infer $\vdash a \to a_1 \to ...\ a_n \to b$, and his amalgamation of identical subcomponents, from $\vdash a \to a \to a_1 \to ...\ a_n \to b$ infer $\vdash a \to a_1 \to ...\ a_n \to b$, as they are formulated by \citet[\S2.6]{landini2012frege}. 
Moreover, \mp~yields Frege's generalized modus ponens, from $\vdash a$ and $\vdash a_1 \to ...\ a_n \to a \to b$ infer $\vdash a_1 \to ...\ a_n \to b$; 
with \hypsyll~we recover Frege's generalized hypothetical syllogism, from $\vdash a_1 \to ...\ a_n \to a \to b$ and $\vdash b_1 \to ...\ b_m \to a$ infer $\vdash a_1 \to ...\ a_n \to (b_1 \to ...\ b_m \to b)$;
we may use \dilemma~to yield any instance of Frege's generalized contraposition, such as from $\vdash a_1 \to ...\ a_n \to a \to b$ infer $\vdash a_1 \to ...\ a_n \to \neg b \to \neg a$;
using \dilemma~we deduce Frege's generalized dilemma, from $\vdash a \to a_1 \to ...\ a_n \to b$ and $\vdash \neg a \to a_1 \to ...\ a_n \to b$ infer $\vdash a_1 \to ...\ a_n \to b$; 
with \ug~we may restore Frege's concavity introduction, from $\vdash a_1 \to ...\ a_n \to f(x)$ infer $\vdash a_1 \to ...\ a_n \to \forall x f(x)$. 
It remains to be seen how the rules of amalgamation of horizontals used by Frege may be recovered. Though rather simple, this result deserves a lemma due to its relevance:

\begin{lemma}[Amalgamation] \label{lem:inner-amal}
	The following rules of inference are admissible:
	\begin{enumerate}
		\item $\Gamma \vdash g(--a) \Longleftrightarrow \Gamma \vdash g(-a)$
		\item $\Gamma \vdash g(-(\neg -a)) \Longleftrightarrow \Gamma \vdash g(\neg a)$
		\item $\Gamma \vdash g(- ((-a) \to (- b))) \Longleftrightarrow \Gamma \vdash g(a \to b)$
		\item $\Gamma \vdash g(- \forall x (- a(x))) \Longleftrightarrow \Gamma \vdash g(\forall x f(x))$
		\item $\Gamma \vdash g(-(a = b)) \Longleftrightarrow \Gamma \vdash g(a = b)$
	\end{enumerate}
\end{lemma}
\begin{proof}
	By {\mp} on {\iiia} and {\ahor--\aid} to treat cases (1)--(5), respectively. 
\end{proof}

It is also worth calling attention to Frege's remark in \gga~\S5 that his turnstile sign is composed of his vertical- and horizontal-stokes. Lemma~\ref{lem:inner-amal} makes it clear that {\GGA} accommodates this view in that the horizontal is contained in the provability turnstile. In fact, this may be expressed as an admissible rule of inference as well: 

\begin{lemma}[Turnstile] \label{lem:turnstile}
	$\Gamma \vdash a \Longleftrightarrow \Gamma \vdash -a$.
\end{lemma}
\begin{proof}
	In the forward direction, by \mp~on \ia, we have $\Gamma \vdash \top \to a$. 
	By Lemma~\ref{lem:inner-amal}~(3) we may unpack the implication and obtain $\Gamma \vdash -(-\top \to -a)$, introduce an additional horizontal by \ref{lem:inner-amal}~(1) to yield $\Gamma \vdash -(-\top \to -(-a))$, and then pack it again into $\Gamma \vdash \top \to -a$. We then apply \mp~and prove $\top$ from \iiie~and~\ug. 
	The backward direction is analogous, except that the application of \ref{lem:inner-amal}~(1) is used to remove and not introduce the extra horizontal. 
\end{proof}

There is another reason why Lemma~\ref{lem:inner-amal} is important. 
Though we show in~\Cref{subsec:completeness} that it is not always the case that $\Gamma \vdash g(a)$ iff $\Gamma \vdash g(-a)$, for all $a \in \term(\sigma)$, 
this lemma reveals that Landini's (\ref{A}) has provable instances whenever the term in question is not atomic:

\begin{corollary} \label{lem:a-cor} 
	If $a$ is not atomic, then $ a \to ((- a) = a)$ is a theorem.
\end{corollary}
\begin{proof}
	By induction on the structure of the term. The base case is vacuously true and the inductive case follows by modus ponens on \ia~and by employing the axioms and rules of amalgamation of horizontals to show that $(- a) = a$ for every connective. 
	Identity and horizontal are immediate from \ahor~and \aid~and the symmetry of identity. 
	For implication, 
	by \aimp~and Lemma~\ref{lem:inner-amal}~(1) we are free to introduce another horizontal stroke on the right-hand side of the identity sign to yield $- (- ((-a) \to (- b))) = a \to b$. Then we apply Lemma~\ref{lem:inner-amal}~(3) once again to transform $- ((-a) \to (- b))$ back into the implication $a \to b$. Negation and the universal quantifier can be treated mutatis mutandis. 
\end{proof}

Lastly, before moving on to the semantics of \GGA, we consider two simple but subtle results about the horizontal that further support our axiomatization. Both will also play a significant role in our completeness proof in~\Cref{sec:completeness}. 
First, the formula $a \to (- a) = \top$ may be regarded as a weaker counterpart of (\ref{A}) in that it states that if $a$ is true then the horizontal~of $a$ is identical to the true. 
Suggestively labeled \eqref{aminus}, we expect this formula to be a theorem of \GGA, given that \citet[p.~352]{landini2022stipulations} sketches a proof of it using Basic Law~IV in the subsystem of \textit{Grundgesetze} without value-ranges. As it turns out, to flesh out his arguments in \GGA, we must rely on the explicit axiomatization of the horizontal, which cannot be easily visualized in Frege's system. 
Second, \citet[\S4]{landini2022stipulations} also shows that (\ref{A}) is provably equivalent to the principle $a \to a = \top$~(A+), a strengthened version of the formula considered above without a horizontal. This too can be shown in \GGA.  
The proofs below rely not just on \iv~but also on \aquant~and, through Lemma~\ref{lem:turnstile}, \ahor~and \aimp. 

\begin{lemma} \label{lem:admissible} 
	The following are theorems:
	\begin{enumerate}
		\item $a \to ((- a) = \top)$ \mytag[aminus]{A-}%
		\item $(a \to (a = \top)) = (a \to ((-a) = a))$ 
	\end{enumerate}
\end{lemma}
\begin{proof}
	We show that $a \to ((- a) = \top)$ is a theorem first. 
	By deduction, we just have to show that $a \vdash (-a) = \top$. If we can introduce a horizontal on the right-hand side of the identity then we can apply \iv. 
	We do this by the transitivity of identity obtained from \iiia~and \iiie~in the usual way. 
	Since $\top$ abbreviates $\forall x (x = x)$, what we have to establish is that $a \vdash \forall x (x = x) = (-(\forall x (x = x)))$. This is a consequence of Corollary~\ref{lem:a-cor}, since, as shown in our proof of Lemma~\ref{lem:turnstile}, $\top$ is a theorem. 
	We now use \mp~on \iv~to derive $a \vdash (-a) = (-\top)$. 
	The rule of negation introduction is admissible, since, as already remarked, a Mendelson-style system is included in {\GGA}. 
	Therefore, it is enough to show that $a, (-a) = \neg \top \vdash \bot$. 
	We instantiate \iiia~with $((-a) = \neg \top) \to (-a \to \neg \top)$ and apply \mp~to obtain $a, (-a) = \neg \top \vdash -a \to \neg \top$. We apply Lemma~\ref{lem:turnstile} to derive $a \vdash -a$. By \mp, $a, (-a) = \neg \top \vdash \bot$ as required.

	We derive (2) by transforming our goal into $\vdash -(a \to (a = \top)) = -(a \to (-a) = a)$ by repeated uses of Lemma~\ref{lem:inner-amal}~(3) and (1). We then further refine it into $-(a \to (a = \top)) = \neg (a \to (-a) = a) \vdash \bot $ by \iv~and negation introduction. By the symmetry of identity, \mp~on \iiia, and deduction, we derive the sequent $-(a \to (a = \top)) = \neg (a \to (-a) = a), \neg (a \to (-a) = a) \vdash -(a \to (a = \top))$. By Lemma~\ref{lem:turnstile}, deduction, and symmetry, this yields~$-(a \to (a = \top)) = \neg (a \to (-a) = a), \neg (a \to (-a) = a), a \vdash \top = a$. Now by \mp~on \iiia~with the identity statement $\top = (-a)$ obtained by \eqref{aminus} and symmetry, we infer 
	$-(a \to (a = \top)) = \neg (a \to (-a) = a), \neg (a \to (-a) = a), a \vdash (-a) = a$. We apply the deduction theorem to immediately arrive at the desired contradiction. 
\end{proof}

\subsection{Semantics} \label{semantics}

We now turn to the semantics of \GGA. To properly model the horizontal, we introduce a class of algebras of boolean truncation, consisting of a truncation function in addition to a domain and interpretation function. 
We say that a triple $\mathcal{M} = \langle D, I, T\rangle$ is a model iff:

\begin{itemize}
	\item $D$ is a pointed set $\langle A, a_1, a_2 \rangle$ with two distinguished elements, which will be called the true and false and respectively denoted by $1_D$ and $0_D$.  We shall not distinguish $D$ and its carrier set to avoid pedantism. 
	We also require functions $not : \bool \rightarrow \bool$, $impl : \bool^2 \rightarrow \bool$, $all : \bool^D \rightarrow \bool$, and $id : D^2 \rightarrow \bool$ for the usual boolean negation and implication operations and for the universal quantification and identity characteristic functions, respectively.

	\item $I$ is an interpretation function in the classical sense, except that since every signature~$\sigma$ only contains function symbols, we just require that for each $n$-ary function symbol $f_n \in \sigma$, we have a total function~$I(f_n) \in D^n \rightarrow D$. 
	
	\item $T : D \rightarrow \{0_D, 1_D\}$ is an idempotent ``boolean-valued'' function such that $T$ has the truth values $0_D$ and $1_D$ as fixed points. In other words, we must have $T(T(e))=T(e)$ for all elements~$e \in D$ and $T(0_D)=0_D$ and $T(1_D)=1_D$. 
\end{itemize}

Intuitively, $T$ is a valuation function that truncates every object in the domain~$D$ to either the true or the false and agrees with the semantics of the logical connectives. Indeed, we may regard $T$ as a homomorphism that preserves the ``boolean structure'' of $D$, given our stipulation that the truth values are fixed points. This means that $T$ is a map that commutes with boolean negation~$not$, implication~$impl$, and universal quantification~$all$. It does not commute with $id$ since it takes arbitrary elements as arguments, and not truth values or properties; but idempotence still applies because $id$ is boolean-valued. 
It follows that, for any truth values $u, v \in \bool$, property $h : D \rightarrow \bool$, and elements $e, d \in D$:
\begin{align*}
	T(not(u))&= not(T(u)) \\						
	T(impl(u,v))&= impl(T(u),T(v)) \\				
	T(all(h))&= all(\fun{e}{T(h(e))}) \\ 		
	T(id(e, d))&= id(e, d). 						
\end{align*}

Some terminology will be useful to categorize the class of models of \GGA. 
We call a model in which $T(e)=1_D$ iff $e=1_D$ for all $e \in D$ a standard model, since this behavior follows Frege's intended interpretation of the horizontal. It should not come as a big surprise to readers that not every model of {\GGA} is standard, given that in \Cref{sec:independence} we established the independence of (A) with respect to the subsystem of {\gga} without value-ranges. We provide in Corollary~\ref{cor:nonthm} our first explicit example of a nonstandard model of \GGA, though in fact our independence proof already loosely points to such a nonstandard model.

We are now prepared to define our interpretation of terms. Since {\GGA} is a term-based language, term evaluation represents our main semantic concept. The idea is to employ $T$ to interpret the~horizontal, thus allowing different interpretations across models, while all other connectives have fixed interpretations. Once again, recall that $\mu, e/x$ denotes the variable assignment that reassigns $e \in D$ to $x$ in $\mu$ but leaves all other variables unchanged:

\begin{definition}[Term evaluation] \label{terminterpretation}
	Given a model $\mathcal{M} = \langle D, I, T \rangle$ and a variable assignment $\mu$, there is a unique term evaluation function $-^{\mu}_\mathcal{M} : \term(\sigma) \rightarrow D$ such that
	\begin{align*}
		{x_i}^{\mu}_\mathcal{M} &= \mu(x_i) \\ 
		f_i(a_1, ..., a_n)^{\mu}_\mathcal{M} &= I(f_i)({a_1}^{\mu}_\mathcal{M}, ..., {a_n}^{\mu}_\mathcal{M}) \\
		(a = b)^{\mu}_\mathcal{M} &= id({a}^{\mu}_\mathcal{M},{b}^{\mu}_\mathcal{M}) \\
		(-a)^{\mu}_\mathcal{M} &= T({a}^{\mu}_\mathcal{M}) \\
		(\neg a)^{\mu}_\mathcal{M} &= not(T({a}^{\mu}_\mathcal{M})) \\
		(a \to b)^{\mu}_\mathcal{M} &= impl(T({a}^{\mu}_\mathcal{M}),T({b}^{\mu}_\mathcal{M})) \\
		(\forall x a)^{\mu}_\mathcal{M} &= all(\fun e T({a}^{\mu, e/x}_\mathcal{M})).
	\end{align*}
\end{definition}

The satisfaction relation $\M \vDash_{\mu} a$ is derived from term evaluation. 
Like provability, this is in fact a relation on formulas,
but, again, we omit the apostrophe of a formula and refer to them as terms to simplify our notation and improve readability. This means that, unlike in modern formula-based logical theories, the satisfaction relation does not support the usual inductive definition on the structure of formulas\textemdash  except in the trivial sense, because a formula can only be obtained in one way by embedding a term. 
To say that a term is true means that its truncation denotes the true. In other words, we may write $\M \vDash_{\mu} a$ as a convenient abbreviation for the more primitive fact that ${T(a^{\mu}_{\M})} = 1_D$. More precisely, we stipulate:

\begin{definition}[Satisfaction] \label{def:satisf}
	$\M \vDash_{\mu} a$ iff ${T(a^{\mu}_{\M})} = 1_D$. 
\end{definition}

We define semantic consequence~$\Gamma \vDash a$ in the usual way. 
We just note that this amounts to $\Gamma \vDash a$ iff for all models $\M$, if ${T(b^{\mu}_{\M})} = 1_D$ for all $b \in \Gamma$ then ${T(a^{\mu}_{\M})} = 1_D$. 
Soundness is immediate for most of our axioms schemes and rules of inference, since $T$ is the identity on the true and the logical connectives are interpreted classically: 

\begin{theorem}[Soundness]
	$\Gamma \vdash a \Longrightarrow \Gamma \vDash a$. 
\end{theorem}
\begin{proof}
	By a straightforward induction on the structure of derivation. It suffices to check that the novel axiom schemes~\iv--\aquant~about the horizontal are valid. Let $\M$ be an arbitrary model and $\mu$ any variable assignment. We can immediately verify the following:
	\begin{itemize}
		\item[\iv.] $\M \vDash_{\mu} (\neg ((-a)=(\neg b)))\to((-a) = (-b))$, as 
		either $T(a^{\mu}_\M) = T(b^{\mu}_\M)$ or not. 
		
		\item[\ahor.] $\M \vDash_{\mu} (-a) = (-(-a))$, because we stipulated that $T(a^{\mu}_\M) = T(T(a^{\mu}_\M))$; 
		
		\item[\aneg.] $\M \vDash_{\mu} (\neg a) = (-(\neg -a))$, since $T$ is idempotent and commutes with $not$, so we see that
		$not(T(a^{\mu}_\M)) = not(T(T(a^{\mu}_\M))) = T(not(T(a^{\mu}_\M))) = T(not(T(T(a^{\mu}_\M))));$
		
		\item[\aimp.] $\M \vDash_{\mu} (a \to b) = (-(-a \to -b))$ for similar reasons as the case above; we only remark that we have  $impl(T(a^{\mu}_\M), T(b^{\mu}_\M)) = T(impl(T(T(a^{\mu}_\M)), T(T(b^{\mu}_\M))));$
		
		\item[\aquant.] 
		$\M \vDash_{\mu} (\forall x f(x)) = -(\forall x (-f(x)))$, this is not only because $T$ is a boolean-structure preserving idempotent homomorphism but also due to function extensionality, which together yield $all(\fun e T({f(x)}^{\mu, e/x}_\mathcal{M})) = T(all(\fun e T(T({f(x)}^{\mu, e/x}_\mathcal{M}))));$

		\item[\aid.] $\M \vDash_{\mu} (a = b) = (-(a = b))$ follows from our stipulation on the fixed points of $T$. 
	\end{itemize}
\end{proof}

We can extract three notable nonderivability results from soundness. The first concerns (A) and the second and third have to do with alternative formulas (A+) and (A*), which are shown by \citet[\S4]{landini2022stipulations} to be provably equivalent to (A) in the subsystem of \gga~without value-ranges (we show this for (A+) in Lemma~\ref{lem:admissible}~(2)). 
It is easy to see that the third formula is but a paraphrase of $(-a) = (a = (a = a))$, and therefore represents a formal embodiment of Frege's famous remark in \S10 that the horizontal function $\GGhorizontal \xi$ can be reduced to the function $\xi = (\xi = \xi)$ using only the identity: 

\begin{corollary}\label{cor:nonthm}
	It is not the case for all $a \in \term(\sigma)$ that we have as theorems: 
	\begin{enumerate}
		\item $a \to ((-a) = a) $ (A) 
		\item $a \to (a = \top) $ (A+)
		\item $(-a) = (a = \top)$ (A*) 
	\end{enumerate}
\end{corollary}
\begin{proof}
	Let $a$ be a constant and build a countermodel~$\M_1$ with $D_1=\{0,1,2\}$, where, as the notation suggests, $0$ and $1$ are the distinguished elements false and true, $I_1(a)=2$, and $T_1(0)=0, T_1(1)=1$, but $T_1(2)=1$. We define the functions $not_1$, $impl_1$, $all_1$, and $id_1$ in the obvious way so they agree with $T_1$.
	Then, $\M_1$ falsifies (\ref{A}), for $\M_1 \vDash_{\mu} a \to ((-a)=a)$ iff $T(a^{\mu}_{\M}) \neq 1 \text{ or } T(a^{\mu}_{\M}) = a^{\mu}_{\M}$. Yet, as it turns out, $T_1(2)=1$ and $T(2) \neq 2$. 
	The other two formulas are invalidated for analogous reasons. 
\end{proof}

The countermodel described in the above corollary is an example of a nonstandard model. In~fact, this is precisely the kind built in~\Cref{sec:independence} in our independence proof.  We can now show that the following rule of inference is not admissible: from $\Gamma \vdash g(a)$ infer $\Gamma \vdash g(-a)$. If $g$ is a function symbol of arity one, for example, just extend the interpretation function in $\M_1$ so that $I_1(g)$ is a function mapping $2$ to $1$ but not $1$ to $1$. Because of the admissibility of the rules amalgamation of horizontals for logical connectives from Lemma~\ref{lem:inner-amal}, however, the above rule is admissible in special cases when $a$ is not atomic. 

\subsection{The completeness proof for \GGA} \label{subsec:completeness}

In the remainder of this section, we give a Henkin-style completeness proof for \GGA.  
To the best of my knowledge, this is the first completeness result for a quantified fragment of \gga~that gives full consideration to the horizontal. 
As $\Gamma \not\vdash a$ implies that $\Gamma \cup \{\neg a \}$ is consistent, 
we just need to show that $\Gamma \cup \{\neg a \}$ has a model. 
Curiously, due to the idiosyncrasies of Frege's logical views on objecthood, our proof of the model existence lemma is rather unique in that the elements of the domain include truth values themselves. 
Moreover, the nonprimitive status of the satisfaction relation simplifies our proof:  
the standard claim that membership in a Henkin set $\Gamma^*$ amounts to truth in the canonical model constructed from $\Gamma^*$ may be proven directly, avoiding the usual argument by induction on formulas. 
In the proof below our use of a term-based system is thus most evident. 
This is also the only place where we rely on our assumption in~\Cref{language} that our signature~$\sigma$ is countable. 

\begin{lemma}[Model existence] 
	If $\Gamma \subseteq \form(\sigma)$ is consistent, then $\Gamma$ has a model. 
\end{lemma}

\begin{proof}
We extend $\sigma$ into a new signature~$\sigma^*$ with countably many fresh constants $c_0, c_1, ...$ 
We then expand $\Gamma$ into a consistent Henkin set $\Gamma^*$ in the usual way. We define an increasing infinite sequence of sets $\Gamma_0 \subseteq \Gamma_1 \subseteq ...$ by induction going through all elements of $\form(\sigma^*)$. 
The construction is the familiar one, except that our focus is on terms. It is easy to see using the standard arguments that this indeed results in a Henkin set. This only relies on the basic first-order fragment \ia--\iia, \hypsyll--\fug, \mp--\ug~without the horizontal and identity. 
Next, we build the canonical model~$\mathcal{M} = \langle D, I, T\rangle$ using the Henkin set $\Gamma^*$:

\begin{itemize}
	\item 
	We start by considering, just as in the ordinary treatment of first-order logic with identity, the set~$E$ of all closed terms of $\Gamma^*$ identified by provable equalities. That is, if $a = b \in \Gamma^*$, we pick the lowest term in the enumeration, unless either is $\top$, in which case we always pick $\top$, since we want it to be our distinguished element $1_D$. 
	However, since $\Gamma^*$ is consistent, clearly $\bot \notin E$. Yet, we need a false~$0_D$ element as well. To address this issue, we define our domain as $D = E \cup \{ \bot\}$ and $\bot$ as $0_D$.

	We let $not$, $impl$, $all$, $id$ be their respective syntactic counterparts $\fun u \neg u$, $\lambda u v. u \to v$, $\fun h \forall x \hat{h}$, $\lambda e d. e = d$, where $\hat{h} := \top$ iff $h(e) = 1_D$ for all $e \in D$ and $\hat{h} := \bot$ otherwise. To prove that these are indeed functions with the required codomain we use \aneg--\aquant, often in the rule forms from Lemma~\ref{lem:inner-amal}~(2--5). 
	As an example, here we will only show that for all $e, d \in D$, $e = d \in \{\bot, \top\}$. 
	Since $\Gamma^*$ is Henkin and thus maximally consistent, either $(e = d) \in \Gamma^*$ or $\neg (e = d) \in \Gamma^*$. 
	If the former, by {\mp}, \iiia, Corollary~\ref{lem:a-cor}, and Lemma~\ref{lem:admissible}~(2), $\Gamma^* \vdash (e = d) = \top$. By closure under derivability, $(e = d) = \top \in \Gamma^*$; By our definition of $D$, this means that $e = d$ and $\top$ are identical elements. 
	Now, if $\neg (e = d) \in \Gamma^*$, we show that $\Gamma^* \vdash (\neg (e = d)) \to ((e = d) = \bot)$. By deduction, we use Lemma~\ref{lem:inner-amal}~(5) to introduce a horizontal on the left-hand side of the outermost identity sign in the consequent, so that it suffices to show that $\Gamma^*, \neg (e = d) \vdash (- (e = d)) = \bot$. 
	As $\bot := \neg \top$, by \mp~on \iv~via contraposition, this refines our goal to deriving the sequent $\Gamma^*, \neg (e = d), (e = d) = \top \vdash \bot$, which can be inferred from \iiia.

	\item We define $I$ in the usual way. Given any $n$-ary $f_i \in \sigma^*$, we let $I(f_i)(a_1,...,a_n)$ be the term $a \in D$ such that $a = f_i(a_1,...,a_n) \in \Gamma^*$, for $a_1, ..., a_n \in D$.

	\item We define $T$ using the horizontal by stipulating that $T(a):=b$, where $b \in D$ is the term such that $b = -a \in \Gamma^*$. That $T$ is indeed idempotent follows from \ahor~plus the fact that identity is an equivalence relation via \iiia~and \iiie. It is easy to see that $T$ has the syntactic truth values~$\top$ and $\bot$ as fixed points. 
	Corollary~\ref{lem:a-cor} and Lemma~\ref{lem:admissible}~(2) yield $\vdash -\top = \top$ by {\iiia} (or as in our proof of Lemma~\ref{lem:admissible}~(1)).  
	By \mp~on \iiia, we may derive $\vdash (\neg -\top) = \bot$ given that $\bot := \neg \top$. 
	We instantiate \aneg~with $\vdash \neg (-\top) = - \neg - (-\top)$. 
	By repeated applications of \iiia~with $\vdash -\top = \top$, we derive $\vdash \neg -\top = - \bot$. We stress that \iv~ensures that the codomain of $T$ is $\{\bot, \top\}$ as we expect. 
\end{itemize}

This shows that~$\mathcal{M}$ is indeed a model. To complete our proof, it is enough to show that, for any possibly open $\sigma^*$-term $a(x_1, ..., x_n)$ depending on the variables $x_1, ..., x_n$, we have that $a(b_1, ..., b_n) \in \Gamma^*$ iff $\M \vDash_{b_1/x_1, ..., b_n/x_n} a(x_1, ..., x_n)$. 
Unfolding Definition~\ref{def:satisf}, we see that our goal amounts to showing that $a(b_1, ..., b_n) \in \Gamma^*$ iff $(-a(b_1, ..., b_n)) = \top \in \Gamma^*$. 
The forward direction follows by \mp~on (A-) from Lemma~\ref{lem:admissible}~(1). 
In the backward direction, we deduce $\top \to -a(b_1, ..., b_n) \in \Gamma^*$ from the symmetry of identity and~\iiia. Since $\top \in \Gamma^*$, by \mp, $-a(b_1, ..., b_n) \in \Gamma^*$. 
We apply Lemma~\ref{lem:turnstile} to obtain $a(b_1, ..., b_n) \in \Gamma^*$. 
\end{proof}

This concludes the proof of our completeness theorem:

\begin{theorem}[Completeness]
	$\Gamma \vDash a$ implies $\Gamma \vdash a$. 
\end{theorem}

It may be helpful to pinpoint how each of the axiom schemes {\iv---\aid} that govern the horizontal contribute to the completeness proof for a better appreciation of this result. Each is used in various ways to keep our construction of the canonical model in check.  {\ahor}, {\aneg}, {\aquant} constrain $T$ to ensure that it is idempotent and has the truth values as fixed points; {\iv} makes sure that truncated terms only have $\top$ and $\bot$ as values. {\iv} and {\aneg--\aquant} are used to show that $not$, $impl$, $all$ and $id$ are their own syntactic counterparts in the formalism. For example, without \aid~our argument that $e = d \in \{\bot, \top\}$ would not go through. {\aimp} is used to collapse horizontals on the right-hand side of the turnstile. This makes sure that provability in the canonical model agrees with our satisfiability relation.

It may also be instructive to take a moment to reflect on why this proof would fail were we to restrict our class of algebras to only standard models.  
We would need to extend $\Gamma$ so that $(-a) \in \Gamma^*$ iff $a = \top \in \Gamma^*$, since in a standard model~$\M$, $\M \vDash_\mu a $ iff $a_\M^\mu = 1_D$. 
But such an extension is not possible because $\Gamma \cup \{a, \neg (a = \top) \}$ is consistent. In the worst case scenario, where $a, \neg (a = \top) \in \Gamma_0$, we could never consistently add $a = \top$ to it. To make the proof go through an additional principle such as (\ref{A}) is needed 
to ensure that $a, \neg (a = \top) \vdash \bot$. All of this motivates the following corollary of our completeness proof: 

\begin{corollary}
	\GGA+(A) is complete with respect to standard models. 
\end{corollary}

Finally, we stress that the completeness result for the propositional fragment of \gga~obtained by~\citet[\S3.5]{cadetpanza2015logical} is consistent with the independence of (\ref{A}) and our completeness result here because identity needs not be interpreted as a biconditional when not flanked by the horizontal on both sides. 

\section{Concluding remarks} \label{conclusion}

It is now time to take stock. 
We have established that the mainstream view maintained by~\cite{parsons1987consistency}, \cite{heck1996consistency}, \cite{burgess1998consistent} and many others that the logic underlying {\gga} is classical is mistaken. 
More precisely, Frege's first-order logical system, as encapsulated in system~\GGA, is complete with respect to a nonclassical semantics. \GGA~is best classified as an extended logic in the sense of \citet[pp.~4--5]{haack1978philosophy}, that is, an extension of classical logic that adds new logical
connectives to its language not definable otherwise, together with new axioms or rules of inference to
govern them. 
The extraclassical nature of \GGA~is derived from the presence of the horizontal. 
Indeed, Corollary~\ref{cor:nonthm} already makes it clear that the horizontal plays an indispensable role in the formalism. 
Frege's suggestion in~\S10 that the horizontal is reducible to identity has no formal attainment. Similarly, our results indicate the claim forcefully made by~\citet[p.~315]{dummett1973fregelanguage} the horizontal is ``wholly superfluous'' as a primitive symbol in Frege's system is plainly wrong. 
I hope to have shown that the received view on the logical dispensability of the horizontal is false. 

One might perhaps object that {\GGA} is not what Frege himself would have accepted as a codification of the first-order fragment of his system in \gga. Surely, there is no reason to think Frege would be willing to isolate a first-order fragment of his logic, given his universalist conception of logical laws~\citep{goldfarb2001frege}. But that is beside the~point. It is enough for our purposes that {\GGA} adapts Frege's first-order fragment without value-ranges to a modern setting we can study metamathematically. 
Granted, as I stressed in \Cref{systemgga1}, I highly doubt that any theory built on top of modern logic can provide an unquestionably accurate rendition of the logical system of \gga, just as no translation can stay completely true to the original. 
But this is not an objection to {\GGA} itself. 
{\GGA} only aims to be an accurate modernization of Frege's first-order logic. We saw in \Cref{systemgga1} that versions of Frege's first-order axioms and rules can all be accommodated in \GGA. My claim that the first-order logical system of {\gga} is nonclassical only hinges on the assumption that these versions translate Frege's original axiom and rules to some extent.

In contrast to \GGA, the stronger system \GGA+(A) considered in the previous section is essentially classical first-order logic in the following sense. The horizontal in this system is indeed semantically irrelevant since the truncation function $T(x)$ is reduced to the boolean-valued function $id(x,1_D)$ following Frege's intended interpretation. 
In other words, we may rewrite a standard model simply as $\langle D, I \rangle$ just as in the classical logic, except, of course, for the fact that the domain is a pointed set with two distinguished elements. (A*) holds and the horizontal is definable via identity, making it logically superfluous in \GGA+(A).\footnote{
	This does not contradict \citepos{bellucci2023frege} thesis that the horizontal is \textit{philosophically} indispensable (due to questions of sense) assuming its intended interpretation.
 }
Frege himself was intending to work within classical semantics. This is clear from how (\ref{A}) is valid under his intended interpretation of the horizontal. But for some reason, Frege hesitated to add enough axioms and rules to his formalism to capture this semantics.

This leads us back to Landini's claim about (\ref{A}) discussed at the beginning of \Cref{sec:independence}. 
\citet[\S4]{landini2022stipulations} argues that ``Frege held (A) and likely thought it is a theorem'' by noting that Frege writes in \S10 that $\xi = (\xi = \xi)$ has the same value as $\GGhorizontal \xi$ for every argument and moreover that the provably equivalent (A*) is accepted in the stipulations at \S10 identifying the true and the false with their respective value-ranges. Landini concludes from this ``that Frege must have imagined that he has (A*)'' and that its unprovability was an oversight. 
I am afraid this conclusion conflates the intended interpretation with the formal system. In \S10, Frege informally acknowledges the truth of (A*) under the intended interpretation, but this does not mean that he regarded it as a theorem. Frege carefully distinguishes truth in the intended interpretation and provability in his system. \citet[pp.~31--32]{heck2012reading} calls attention to a footnote to \gga~\S114, where Frege says that a proposition needed to show that there are infinitely many numbers ``is, it seems, unprovable'' but does not add that it is false. As Heck suggests, it is likely that Frege thought of it as true but unprovable and that Frege was prepared to acknowledge that his system was negation incomplete. To sum up, I believe Frege deliberately left out (\ref{A}) from the system of {\gga} simply because it is not needed for the formalization of arithmetic he needed to vindicate his logicism.

This is not to say that (\ref{A}) is not useful to prove certain facts about numbers. 
As Landini points out, having this formula as an additional axiom allows us to have $\neg 2$ as a theorem. But such a theorem would in no way contribute to a formalization of actual arithmetic, 
no more than the theorem $3 \in 17$ about von Neumann ordinals in $\textsf{ZFC}$ can, for example.\footnote{
	\citet[\S4]{landini2022stipulations} argues that (\ref{A}) does not ``overdetermine'' the true as a logical object, in 
	the sense that there are no ``rival logical options available concerning candidate
	axioms which could, for all we know, equally well consistently govern a logical object'' (p.~353). 
	I am not sure about the true, but, as I see it, (\ref{A}) does ``overdetermine'' the natural numbers as logical objects due to the introduction of junk theorems like $\neg 2$. In fact, a very reasonable rival option to (\ref{A}) is not to postulate it at all for a more finely grained logical account of numbers. }  
My main takeaway here is that, due to the absence of (\ref{A}), the system that Frege thought was sufficient to vindicate his logicism ended up being too weak to capture his intended interpretation, and its first-order fragment, encoded by \GGA, turned out to be extraclassical as a result.

Let me conclude by addressing another possible objection. Since {\GGA+(A)} yields classical semantics, one may argue that missing something required for a system to function classically is rather different from claiming that the system is extraclassical. 
I see the extraclassical nature of {\GGA} as comparable to that of a modal normal system such as $\mathsf{T}$. Once $\phi \to \Box \phi$ is added to the system, $\mathsf{T}$ reverts back to what is logically equivalent to classical logic. 
But this does not mean that $\mathsf{T}$ itself is not extraclassical. \GGA~is extraclassical in a similar way. In fact, the term ``extraclassical'' is coined by \citet[ch.~1, \S1.1]{burgess2009philosophical} to refer to modal logics.

\begin{acknowledgments}
	I wish to thank Ansten Klev, Marco Panza, Matthias Schirn, Johan van Benthem, and three anonymous referees for comments on early drafts. 
	This paper was inspired by conversations with 
	Alessandro Bandeira Duarte and Gregory Landini in the summer of 2024. 
	I am specially indebted to Alessandro for bringing Axiom (A) to my attention. 
	The material this paper is based on was presented at the International Academy for Philosophy of Science Conference (AIPS 2025) in Hangzhou on October 28, 2025. I am thankful for questions by Bo Chen, Alan H\'ajek, and Johan van Benthem in the audience. 
\end{acknowledgments}

\bibliographystyle{plainnat}
{\linespread{0.7}\selectfont\bibliography{../../../mybib}}

\address

\end{document}